\documentclass{amsart}

\usepackage{mathrsfs}

\usepackage{amsfonts,amssymb,amsmath,amsthm}

\usepackage{color}

\usepackage{enumerate}

\def\NN{{\mathbb N}}

\def\rho{\varrho}

\def\phi{\varphi}

\newcommand{\be}[1]{\begin{equation}\label{#1}}

\newcommand{\ee}{\end{equation}}

\def\cal{\mathscr}

\def\NN{{\mathbb N}}

\newtheorem*{Thm}{Theorem}
\newtheorem*{cor}{Corollary}

\makeatletter

\@addtoreset{equation}{section}
\makeatother

\makeatletter
\def\blfootnote{\xdef\@thefnmark{}\@footnotetext}
\makeatother

\title[Powered numbers]
{On the sum of two powered numbers} 

\author{J. Brüdern and O. Robert}

\begin{document}


\maketitle


 Many of the historic diophantine problems concern additive properties of powers of natural numbers. Fermat's Last Theorem, the Catalan conjecture and problems of Waring's type are among the most familiar representatives in this class of questions. In his beautiful survey article Mazur \cite{M} suggests to consider rounded versions of these problems. Rather than  insisting  to work with perfect powers, one is then led to analyse the additive properties of natural numbers $m$ whose kernel $k(m)$, defined as the largest squarefree divisor of $m$, is small in terms of $m$. To make this precise, fix a
real number $\theta\in(0,1]$ and define the set 
\[
\cal{A}(\theta)=\{m\in \NN\colon  k(m)\le m^{\theta}\}.
\] 
For $\theta=1/l$ with $l\in \mathbb N$ this set constitutes the powered numbers, in Mazur's terminology. The motivation is that
$ i(m) = \log m/\log k(m)$
measures multiplicities in the prime factorisation of $m$,  and for all $l\in\mathbb N$ one has $i(m^l)\ge l$ so that the $l$-th powers form a subset of $\cal{A}(1/l)$, and elements of $\cal{A}(1/l)$ should then substitute for $l$-th powers in Mazur's rounded diophantine problems.    
As demonstrated in \cite{BR}, the number $S_\theta(x)$ of elements in $\cal{A}(\theta)$ that do not exceed $x$ satisfies the inequalities
\be{simple} x^\theta \ll S_\theta(x) \ll x^{\theta+\varepsilon} \ee
whenever $\varepsilon$ is a given positive real number and $x$ is large in terms of $\varepsilon$. Hence, if one replaces $l$-th powers by members of 
$\cal{A}(1/l)$, then the new set is roughly of the same density.

In this short communication, we discuss the binary linear equation $u+v=n$ with the variables $u,v$ restricted to $\cal{A}(\theta)$. Thus, we 
are interested in the sumset 
\[
\cal{A}(\theta)+\cal{A}(\theta):=\{u+v\colon u,v\in \cal{A}(\theta)    \},
\]
and in particular, we wish to determine the values of $\theta$ where this sumset 
contains all but finitely many natural numbers. Our main result comes close to the rounded version of representations by sums of two squares.

\begin{Thm} 
Let $C=2\sqrt[4]{27}$. Then, for all natural numbers $n$ with  $n\ge 4$ there are natural numbers $m_j$ with $m_j\ge 2$ $(j=1,2)$ and
\be{A}  n=m_1+m_2, \quad  k(m_j) \le C\sqrt{m_j} \quad (j=1,2).\ee
\end{Thm}

This conclusion has an immediate consequence for our initial question.

\begin{cor} Let $\theta>\frac12$. Then all large natural numbers are the sum of two numbers in $\cal{A}(\theta)$.
\end{cor}

Since there are certainly no more than $S_\theta(x)^2$ natural numbers below $x$ in the sumset $\cal{A}(\theta)+\cal{A}(\theta)$, it follows from \eqref{simple} that for all $\theta<\frac12$, the sumset $\cal{A}(\theta)+\cal{A}(\theta)$ has density $0$ in the natural numbers. Therefore, the conclusion in the corollary is certainly false for $\theta<\frac12$.

Our theorem fails to settle the case $\theta=\frac12$. However, the set $\cal{A}(\frac12)$ is somewhat denser than the set of squares, and it is very likely that the condition on $\theta$ in the corollary can be relaxed to $\theta\ge \frac12$. It is worthwhile to consider the vincinity of $\frac12$ on a finer scale. In our companion paper \cite{BR} we established asymptotic formulae for $S_\theta(x)$ and related counting functions. As a consequence of \cite[(1.15)]{BR}, whenever $\varepsilon>0$ and $x$ is large in terms of $\varepsilon$, one has
$$ S_\frac12 (x) \gg \sqrt{x} \exp\big((2-\varepsilon)\sqrt{\log x/ \log\log x}\big). $$
Further, as a special case of \cite[Theorem 3]{BR}, for any real number $\gamma$, one finds
\be{log} \#\{m\le x \colon k(m) \le \sqrt{m} (\log m)^\gamma\}
= (\log x)^\gamma S_\frac12(x)\Big(1+ O\Big(\sqrt{\frac{\log\log x}{\log x}}\Big)\Big). \ee
In light of these estimates, it seems safe to conjecture that for any given real number $\gamma$, a sufficiently large natural number $n$ has a representation $n=m_1+m_2$ with
\be{logss} k(m_j) \le \sqrt{m_j} (\log m_j)^\gamma\quad (j=1,2)\ee
We expect the expansion \eqref{log} to persist if the role of the logarithm in the factors $\log m$ and $\log x$ is played by a function that grows somewhat faster, and one may then hope to be able to tighten the condition \eqref{logss} to 
$$ k(m_j) \le \sqrt{m_j}\exp\big(- \sqrt{\log m_j/\log \log m_j}) $$
or thereabout, and still conclude that $n=m_1+m_2$ has solutions with this condition in place, at least when $n$ is large.
 
\smallskip

The proof of the theorem is entirely elementary. 
We write $A= \frac1{2} \sqrt[4]{27}$. Let $a$ and $b$ be the integers defined by
\be{Defa} 2^a < \sqrt{n}/A \le 2^{a+1}, \quad 3^b < A\sqrt{n} \le 3^{b+1}. \ee
Then $a\ge 1$ and $b\ge 1$ for $n\ge 7$. Next, define $U$ and $V$ by
\be{1} U = [n/2^a] -1, \quad n=2^aU+V . \ee
Then $2^a U = n- \xi 2^a$ for some real number $\xi\in[1,2)$, and so
$ 2^a \le V < 2^{a+1}. $
Since $2$ and $3$ are coprime, there are integers $W$, $w$ with 
\be{2} V=-2^aW+ 3^b w.\ee
The theory of linear diophantine equations shows that in this last equation it can be arranged that $1\le w \le 2^a$. With this choice for $w$ we put
\be{defuv} u=2^a(U-W), \quad v=3^b w. \ee
By \eqref{1} and \eqref{2}, we have $n=u+v$. Further, the constraints on $w$ and the inequalities \eqref{Defa} ensure that
\be{4} 3^b \le v \le 2^a 3^b  <n . \ee
By \eqref{4}, we first see that $u\ge 1$, and then by \eqref{defuv} that $u\ge 2^a$. It only remains  
to estimate the kernels of $u$ and $v$. By construction, we have
$$ k(v) \le 3k(w) \le 3w  = \sqrt{v} \, \frac{\sqrt w}{3^{b/2-1}}. $$
By \eqref{Defa} and the range for $w$, we see 
$$ 3^{2-b} {w}< \frac{2^a}{ 3^{b-2}} < \frac{\sqrt n/A}{A\sqrt n /27} \le \frac{27}{A^2}.   $$
This yields
$$ k(v) < \frac{\sqrt{27}}{A} \sqrt v  . $$

\medskip

The treatment of $u$ is similar. We use $u=2^a(U-W) = n-v$. Then
$$ k(u) \le 2k(U-W) \le 2(U -W),$$ and hence  
$$ k(u) \le 2^{1-a}u = \sqrt u \,\sqrt{n-v}\, 2^{1-a}.$$
But $2^{1-a} \le 4A/\sqrt n$, and then
$$ k(u) \le 4A \sqrt u. $$ 
Now choose $m_1=u$ and $m_2=v$ to  confirm \eqref{A} whenever $n\ge 7$. For  $4\le n\le 6$ one may use $m_1=2$ and $m_2= n-2$.

\vspace{2ex}\noindent
 J\"org Br\"udern \\
 Universit\"at G\"ottingen\\
Mathematisches Institut\\
Bunsenstrasse 3--5\\
D 37073 G\"ottingen \\
Germany\\
jbruede@gwdg.de\\
[2ex]
Olivier Robert\\
 Universit\'e de Lyon,  Universit\'e de Saint-\'Etienne\\
Institut Camille Jordan CNRS UMR 5208\\
25, rue Dr R\'emy Annino\\
F-42000, Saint-\'Etienne\\
 France\\
olivier.robert@univ-st-etienne.fr

\begin{thebibliography}{9}
\bibitem{BR} J. Br\"udern and O. Robert. On the distribution of powered numbers. Forthcoming.
\bibitem{M} B. Mazur. Questions about powers of numbers. Notices Amer. Math. Soc. 47 (2000), no. 2, 195--202.
\end{thebibliography}
\end{document}